\documentclass[preprint,review,11pt]{elsarticle}

\usepackage{natbib}
 \bibpunct[, ]{(}{)}{,}{a}{}{,}%
\usepackage{comment}
\usepackage{subfigure}
\usepackage{floatrow}
\usepackage{lscape}
\usepackage{newunicodechar}
\newunicodechar{✓}{\checkmark}
\usepackage{enumerate}
\usepackage{mathrsfs}
\usepackage[justification=centering]{caption}
\usepackage{algorithm,setspace}
\usepackage[noend]{algpseudocode}
\usepackage{epstopdf}
\usepackage{hyperref}
\usepackage{amsmath}
\usepackage{amssymb}
\usepackage{mathtools}
\usepackage{cases}

\makeatletter
\def\@author#1{\g@addto@macro\elsauthors{\normalsize%
    \def\baselinestretch{1}%
    \upshape\authorsep#1\unskip\textsuperscript{%
      \ifx\@fnmark\@empty\else\unskip\sep\@fnmark\let\sep=,\fi
      \ifx\@corref\@empty\else\unskip\sep\@corref\let\sep=,\fi
      }%
    \def\authorsep{\unskip,\space}%
    \global\let\@fnmark\@empty
    \global\let\@corref\@empty  
    \global\let\sep\@empty}%
    \@eadauthor={#1}
}
\makeatother

\newtheorem{theorem}{Theorem}

\journal{TBD}

\begin{document}

\begin{frontmatter}



\title{Dynamics of Equity, Efficiency, and Efficacy in Home Health Care with Patient and Caregiver Preferences}
\author{Mohammad Firouz\corref{cor1}\fnref{fn1}}
\cortext[cor1]{Corresponding Author; E-mail: mfirouz@uab.edu; Phone: (205)934-8830}
\author{Linda Li\fnref{fn2}}
\author{Daizy Ahmed\fnref{fn3}}
\author{Abdulaziz Ahmed\fnref{fn4}}

\address[fn1]{Department of Management, Information Systems, and Quantitative Methods, University of Alabama at Birmingham, Birmingham, AL 35294.}
\address[fn2]{College of Business, Missouri State University, Springfield, MO 65897.}
\address[fn3]{Department of Mathematics and Statistics, Sri Mayapur International School, Sri Mayapur Dhama, Nadia District, West Bengal 741313.}
\address[fn4]{Department of Health Services Administration, University of Alabama at Birmingham, Birmingham, AL 35233.}

\begin{abstract}

\noindent In this paper we discuss \dots
\end{abstract}


\end{frontmatter}


\pagebreak

\section{Introduction}
\label{sec:Intro}

Over the last few decades, Home Healthcare (HHC) services have significantly increased. HHC offers care at the patient's location by a caregiving team and is presented as an attractive alternative to traditional hospitalization. HHC offers intensive support and care such as daily living services, complex care services, and social and psychological assistance. The Center for Disease Control and Prevention (CDC) says that 1 in 4 adults has a disability that impacts major life activities. The possibility of receiving care within one's own home compared to a hospital environment has motivated patients to opt for such services more and more. At the same time, decreasing overhead costs and potentially releasing some of the resources that would have been otherwise occupied by hospitalization of the patients, has made HHC an attractive option for the care providers too. Being traditionally dominated by the non-profit organizations, HHC has now attracted a wide attention from the for-profit firms as well \citep{rosenau2001comparative}. 

The main resource involved in HHC is the team of caregivers which can be volunteers or paid workers. Today, more than 1 in 5 Americans (21.3\% of the population) are caregivers, having provided care to an adult or child with special needs at some time in their past 12 months. This totals to an estimated 53.0 million adults in the United States. On average, a family caregiver will spend over 24 hours each week providing care to someone, although many report spending over 40 hours per week on caregiving duties. This imbalance in the workload on the caregivers' side has led to numerous burn-out situations, where the caregiver is overwhelmed with managing the responsibilities of the personal life versus the care-related ones \citep{NACCareGiving}. 

Every care situation is unique. CDC tracks six disability types:
(1) mobility (serious difficulty walking or climbing stairs),
(2) cognition (serious difficulty concentrating, remembering or making decisions),
(3) hearing (serious difficulty hearing),
(4) vision (serious difficulty seeing),
(5) independent living (difficulty doing errands alone),
(6) self-care (difficulty dressing or bathing). The level of need of each patient for each one of these care services generally is dependent on the level of their age and other health-related measures. While changing over time, the healthcare providers have prior information on the level of care services patients generally need at each stage of their life. The level of expertise or willingness of each caregiver in each one of these services also differs as well, but it can be known generally when a caregiver signs up for the services. The assignment of the right caregiver to the right patient with the objective of increasing patient satisfaction therefore becomes a crucial task. 

There are three main entities in an HHC setting: the firm, the caregivers, and the patients. Often times the interests of each group of entities conflict with the other two. Whether a non-profit or a for-profit HHC, the firm generally seeks efficiency in terms of maximum utilization of the resources and highest level of care provided to the patients. Caregivers, seek a fair utilization across the network in order to increase job satisfaction. Finally, the patients seek maximum satisfaction level and a fair share of the care provided by the firm. Therefore, looking at HHCs from either one of these lenses by themselves leads to decisions that may neglect the interests of the others. 

A common goal among both for-profit and non-profit HHCs, efficiency plays a crucial role in the profitability of the former and sustainability of the latter \citep{rosenau2001comparative}. Additionally, the increasingly patient-driven decision making in healthcare prohibits waste of resources which could have served increase patient satisfaction. Major reliance on volunteer caregivers for non-profit HHCs and  being the most expensive resource in the for-profit HHCs, make waste of manpower as one of the most undesirable practices in an HHC firm. Therefore, one of the goals of this research is to address waste of manpower with a particular focus on maximizing patient satisfaction rates. 

Health equity has a long-standing presence in healthcare research with major focus to health access. Equity in healthcare refers to the absence of systematic disparity in health access by sufficiently different group of individuals. In the for-profit sector the idea of health equity is present in the form of regulations and mandates which often conflicts with efficiency goals. On the other hand, in the non-profit sector, scarcity of the resources and manpower makes the objective of meeting the care demands irrelevant. What remains crucial in non-profit settings, however, is the equitable distribution of the available resources to the patients which is a distinguishing aspect of non-profit HHCs compared to their for-profit counterparts \citep{ellenbecker1995profit}. In other words, although equity is a mandate to the for-profit HHCs, it is a goal for non-profit ones. Therefore, a second goal of this research is to provide a structure by which equitable distribution of resources by an HHC organization to the public is addressed whereby each patient is treated in a fair and equitable manner in terms of assigned services. Within the context of for-profit HHCs, this goal serves as the analysis on the price of equity as meeting certain level of equity generally means sacrificing a certain level of profitability. On the other hand, equity is generally in conflict with minimization of manpower waste in the non-profit HHCs. 

Additionally, National Alliance for Caregiving (NAC) states that nearly two-thirds of caregivers are volunteers who work other full-time or part-time jobs. Inattention to the needs of the caregivers in both for-profit and non-profit HHCs has caused many issues such as sacrifice of personal time, imbalanced personal life, as well as physical, emotional, and psychological burn-outs \citep{NAC}. In the for-profit sector, over-utilization of  caregivers has been reported to lead to long-term physical, emotional, and psychological damages, while in the non-profit sector neglecting volunteer needs, especially with regards to their over and under utilization, has been presented as the major cause of failure in volunteer retention by non-profit organizations leading to ineffective management of manpower \citep{clary1992volunteers}. Although, effective utilization by moderating over and under utilization across the caregivers may seem to be conflicting in the short-term with the other two objectives of efficiency and equity, in the long-term it improves employee well-being, which in turn can significantly aid the other two goals. Therefore, the third goal of this research is to address this neglected issue in HHC research by providing an effective utilization of manpower across the network. 

\cite{beamon2008performance} define three pivotal performance measures of equity, efficiency, and efficacy to gauge  humanitarian relief operations. In our setting, we define the following performance measures accordingly: 
\begin{itemize}
	\item Equity: For each patient, deviation from their ``fair'' share of service should be minimized. 
	\item Efficiency: the maximum possible hours of care should be used across the network.
	\item Efficacy: for each caregiver, deviation from their ``fair'' utilization should be minimized.
\end{itemize} 

Specifically, equity in a caregiver network ensures that each patient is achieving the same satisfaction rate as far as possible. Satisfaction rate of each patient is directly proportional to the number of care hours they receive. On the other hand, efficiency ensures that the assignment policy has the lowest possible waste of the total available care hours across the network which in turn ensures that the maximum number of available care hours is assigned to the patients. Finally, efficacy requires that each caregiver receive their fairest possible utilization of time so that the possibility of burn-outs due to over-utilization and demotivation due to under-utilization is minimized across the network, which can potentially lead to an ineffective policy.

\begin{figure}[H]
	\centering
	\includegraphics[width=0.3\textwidth]{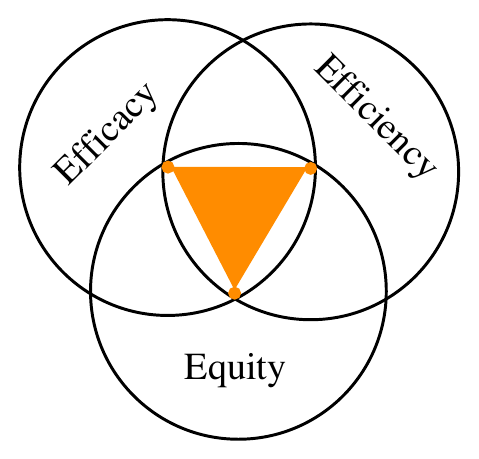}
	\caption{The dynamics of equitable, efficient, and effective policies.}
	\label{fig:Equity_Efficiency_Effectiveness}
\end{figure}

Figure \ref{fig:Equity_Efficiency_Effectiveness} shows the dynamics of the three performance measures of equity, efficiency, and efficacy in HHC settings. The center of each circle represents the perfectness when each measure is considered by itself while the highlighted triangle represents the area in which the HHC firm in our setting may choose to perform in, depending on the flexibility of their policy regarding each one of the three performance measures as well as practical considerations in their network. 

The remainder of this paper is organized as follows: Section \ref{sec:LitReview} reviews the research most related to our paper. Section \ref{sec:Model} introduces our equity-efficiency-efficacy model to the HHC problem. Section \ref{sec:ModAnalysis} describes some of the properties of our problem. Section \ref{sec:NumStudy} highlights our managerial insights related to equity, efficiency, and efficacy and their interactions in HHC. Finally, Section \ref{sec:Conclusions}, concludes the paper with critical insights from the study and future extensions of our work.
\section{Relevant Literature}
\label{sec:LitReview}

The HHC problem has been studied in the literature by a number of researchers from the perspective of three operational level sub-problems: (i) assignment of caregivers to patients; (ii) scheduling of the caregiver visits; and (iii) routing of caregivers across the network of patients.  \cite{begur1997integrated} provide a decision support system for a non-profit caregiver provider in order to address the problem of minimizing waste in the scheduling of the patient visits. \cite{de1998planning} studies a resource allocation problem for HIV patients with uncertain number of patients, considering the stages of health that they go through. \cite{carello2014cardinality} provide a robust formulation model to ensure continuity of care requirement for different classes of patients with variable demands. 
In order to manage the uncertainty arising from cancellation and arrival of new demands, \cite{cappanera2018demand} provide a robust approach to address the joint problem of assignment, scheduling and routing of an HHC organization. 
\cite{restrepo2020home} solve the two-stage stochastic integrated staffing and scheduling problem for an HHC and investigate the benefits of flexibility in taking recourse actions. A number of other researchers provide different methodologies addressing a subset of the three operational level sub-problems \citep{grenouilleau2020new, zheng2021stochastic, lin2021matching, malagodi2021home, grenouilleau2019set, gomes2019modelling, nikzad2021matheuristic}.  

The central theme in all of the above-mentioned papers is that the problem is studied from the context of the HHC firm, neglecting the needs of the patients as well as the caregivers. Therefore, goals of equity of service as well as caregiver retention policies, and their interaction with the efficiency goal are not studied in the literature \citep{rosenau2001comparative}. \cite{carello2018trade} proposes that there are three main stakeholders in an HHC setting: patients, operators, and service providers, and that each party has their own interests that generally conflicts with the others. Patients seek for high quality service, operators seek for fair workloads, and managers seek to reduce costs. They emphasize that models to integrate these three stakeholders in the context of HHC research are of significant interest.

Increasingly patient-focused policies and employee rights regulations derive for-profit HHC firms to look at the goals that traditionally may have been common only in the non-profit settings \cite{frank1994nonprofit}. Particularly the need for understanding the trade-off between efficiency and equity has been emphasized in the literature \citep{savas1978equity, kaplan2002allocating, bertsimas2011price, bertsimas2012efficiency}. A particular stream of research to address this need in the literature has been addressed in the context of humanitarian and disaster relief operations \citep{beamon2008performance, orgut2016achieving, huang2012models, taskin2010inventory}.
Our problem differs from this stream of research as ours deals with the day-to-day operation of the HHC as opposed to the one-time nature of the disaster relief operations. 

Equity in healthcare refers to the absence of systematic disparities in provided health to different groups of individuals \citep{braveman2003defining}. Although the measurement of equity in a quantitative fashion can be defined in various terms \citep{marsh1994equity}, what remains constant is that the goals of efficiency and equity are generally competing objectives and quantifying the equity-efficiency trade-off, provides extremely valuable tools for decision makers \citep{bradley2017operations}. \cite{cleary2010equity} use an optimization model to achieve priorities of HIV care while maintaining equitable treatment opportunities across the public. \cite{wilson2005designing} reflect the equity measure in their objective function to provide equal opportunities for HIV patients. Equal opportunity access to health facilities spread over a geographical area is the central theme of \cite{mccoy2014using}. Similarly, \cite{ares2016column} determine equitable and efficient policies to locate roadside clinics for truck drivers in Africa. They find out that minor losses on the side of efficiency can lead to tremendous gains in equity measure across the network.  \cite{bertsimas2013fairness} use equity and efficiency performance measures to allocate deceased donor kidneys to a national patient list. \cite{su2006recipient} design a kidney allocation policy that takes the preferences of the recipients into account while simultaneously exploring the efficiency-equity trade-off. They find out that a choice-based allocation mechanism leads to outcomes in the middle of the equity-efficiency spectrum. 
 
Majority of the healthcare research that has been done with consideration of equity-efficiency trade-off considers national or international policy implementations. However, our work focuses on the operational level decision making of the HHC firms which aims to provide a tool and set of guidelines for day-to-day decision making. In addition, dominated by mainly volunteer care-givers, an aspect that has been widely neglected in the HHC research, is the mismanagement of the caregivers in terms of utilization. Specifically, it has been emphasized that both under-utilization and over-utilization of the caregivers can lead to ineffective policies, leading to long-term issues faced by the HHC firms. In detail, over-utilization of a portion of the caregivers has proven to lead to problems that create chronic damages to the organization arising from volunteer burn-outs \citep{denton2002job}. On the other hand, it has been also emphasized that under-utilization of the care-givers can actually hurt the organization by decreasing the motivation and the sense of creating a change from the volunteer side \citep{clary1992volunteers}. Therefore, the need for models that allow for the choice of homogeneous utilization of the volunteers is increasingly important. This research is taking a step in this direction while simultaneously addressing the equity-efficiency trade-off within the context of HHC research. 
\section{Model Formulation}
\label{sec:Model}
We consider the problem faced by an HHC firm that wants to assign $n$ caregivers to $m$ patients who require $s$ types of care services. Not all caregivers are fully capable of doing all the services. Specifically, $e_{ik}$ is the capability of expertise of caregiver $i$ in service $k$. Caregiver $i$ could be assigned to a patient if s/he is capable of the assigned service, that is, $e_{ik}=1$. Therefore, $e_{ik}$ is a binary variable. The service hours of the day required by patient $j$ for care service $k$ is $D_{jk}$, which is integer.

The maximum number of hours a caregiver can be utilized is $H_i$. Hence, utilization of each caregiver can be defined as the ratio of the total hours assigned for services and $H_i$. If all caregivers have exactly same level of utilization, then the service requirement is equally assigned to all caregivers. For the non-profit HHC, they do not expect over-utilization of some of the volunteered caregivers or idleness of the others \citep{agostinho2012analysis}. So the smoothness of utilization is one of the objectives for the HHC. We define this by efficacy in our paper. On the other hand, fill-rate of each patient is calculated by the total amount of care hours assigned to them for all the services they need divided by their demanded hours of service. Similar to utilization, if all patients receive exactly same fill-rates, then the perfect equity is achieved and the assignment serves the patients with the most fairness. Of course, the equity must be the HHC's goal. However, pursuing efficacy and equity may sacrifice the efficiency. The three have an obvious trade-off because of limited resources. Efficiency is defined by the patient utility in our paper. We borrow this idea of utility from economic concepts \citep{harsanyi1953cardinal}. The reason is because, for the non-profit HHC, the utility could be the welfare of the patients; and for the profitable HHC, the utility could be the profit generated by providing the care service. Using the idea of utility will make our model work in a general way. 

Combining the three objectives (efficacy, equity and efficiency) discussed above, we build a mix-integer-programming model to maximize the total patient utility while simultaneously achieving equity among the care received by the patients and a fair utilization of the caregivers across the network to minimize under and over utilization. Additionally, the preferences of both the patients (in terms of the total number of caregivers assigned to them) as well as the caregivers (in terms of the total number of patients assigned to them) are to be taken into account in the model's constraints. The budget that the HHC could spend is also considered in the constraints. The HHC of course has limited budget and the cost varies depending on the caregiver's skill set and different services performed.

Table \ref{tab:param} lists the notations along with their definitions that we introduced so far and additional ones throughout the rest of the paper. 
\renewcommand{\arraystretch}{1.2}
\begin{table}[H]
\begin{tabular}{ll}
\hline
\textbf{Parameters} & \textbf{Definitions} \\
\hline
$n$ & number of caregivers. 		\\	
$m$ & number of patients.	\\	
$s$ & number of care services.	\\
$\mathcal{I}$ & set of caregivers, $\mathcal{I}:\{1,\dots,n\}$. 	\\	
$\mathcal{J}$ & set of patients, $\mathcal{J}:\{1,\dots,m\}$.	\\	
$\mathcal{K}$ & set of care services, $\mathcal{K}:\{1,\dots,s\}$.\\
$N_j$ & maximum number of caregivers that can be assigned to patient $j$. 	\\	
$M_i$ & maximum number of patients that can be assigned to caregiver $i$.	\\
$H_i$ & maximum time units caregiver $i$ could be utilized.	\\
$C$ & the total assignment budget of HHC.	\\
$c_{ijk}$ & unit assignment cost of caregiver $i$ to patient $j$ for care service $k$. \\
$e_{ik}$ & 1 if caregiver $i$ can perform care service $k$, 0 otherwise. 	\\
$D_{jk}$ & required time of care service $k$ by patient $j$.	\\
$p_{j}$ & the utility of patient $j$ if s/he receives full service. 	\\
$\beta_{j}$ & fill-rate of patient $j$, $\beta_{j} =\frac{\sum_{i\in \mathcal{I}}\sum_{k\in \mathcal{K}}x_{ijk}}{\sum_{k\in \mathcal{K}}D_{jk}}$.\\
$\theta$ & penalty for deviations from the maximum fill-rate $\beta$.	\\
$u_i$ & utilization rate of caregiver $i$, $u_i =\frac{\sum_{j\in \mathcal{J}}\sum_{k\in \mathcal{K}}x_{ijk}}{H_i}$.\\
$\alpha$ & penalty for deviations from the maximum utilization rate $u$. 	\\
\hline
\textbf{Decision Variables} & \textbf{Definitions} \\
\hline
$x_{ijk}$ & units of caregiver $i$'s time assigned to patient $j$ for care service $k$. \\
$z_{ij}$ & 1 if caregiver $i$ assigned to patient $j$, 0 otherwise. \\
$\beta$ & maximum fill-rate across the patients, $\beta=\max_{j\in \mathcal{J}} \beta_j$. \\
$u$ & maximum utilization rate across the caregivers, $u=\max_{i\in \mathcal{I}} u_i$. \\
\hline
\end{tabular}
\caption{Summary of the mathematical notation.}
\label{tab:param}
\end{table}

Our mathematical formulation is as follows.
 
\begin{align}
Max \hspace{3pt} \sum_{j\in\mathcal{J}} p_{j}\beta_{j} -\sum_{j\in\mathcal{J}} \theta (\beta - \beta_j)- \sum_{i\in\mathcal{I}} \alpha(u-u_i) \label{eq:ObjFnc}
\end{align}\setcounter{equation}{0} 
\begin{subnumcases}{s.t.}
   \sum_{i\in\mathcal{I}}x_{ijk} \leq D_{jk} &$\forall j\in \mathcal{J}, \, \forall k\in \mathcal{K}$ \label{eq:Const1}\\
   \sum_{i\in \mathcal{I}} \sum_{j\in \mathcal{J}}\sum_{k\in \mathcal{K}} c_{ijk} x_{ijk} \leq C & \label{eq:Const3} \\
   x_{ijk} \leq e_{ik}z_{ij}H_{i}  & $\forall i\in \mathcal{I}, \, \forall j\in \mathcal{J}, \, \forall k\in \mathcal{K}$ \label{eq:Const4} \\
   z_{ij} \leq \sum_{k\in\mathcal{K}}x_{ijk}& $\forall i\in \mathcal{I}, \, \forall j\in \mathcal{J}$ \label{eq:Const5} \\
   \sum_{i\in\mathcal{I}}z_{ij} \leq N_j  & $\forall j\in \mathcal{J}$ \label{eq:Const6} \\
   \sum_{j\in\mathcal{J}}z_{ij} \leq M_i  & $\forall i\in \mathcal{I}$ \label{eq:Const7} \\
   \beta_j= \frac{\sum_{i\in \mathcal{I}}\sum_{k\in \mathcal{K}}x_{ijk}}{\sum_{k\in \mathcal{K}}D_{jk}} & $\forall j\in \mathcal{J}$ \label{eq:Const8} \\
   u_i=\frac{\sum_{j\in \mathcal{J}}\sum_{k\in \mathcal{K}}x_{ijk}}{H_i}  & $\forall i\in \mathcal{I}$ \label{eq:Const9}\\
   \beta_j \leq \beta  & $\forall j\in \mathcal{J}$ \label{eq:Const10} \\
   u_i \leq u  & $\forall i\in \mathcal{I}$ \label{eq:Const11}\\
   0\leq u,\beta\leq 1;  &  \notag \\
   x_{ijk}\in \mathbb{N}_0; \; z_{ij} \in \{0,1\}& $\forall i\in \mathcal{I}, \, \forall j\in \mathcal{J}, \, \forall k\in \mathcal{K}$, \label{eq:Const12}
\end{subnumcases}
where $\mathbb{N}_0$ is the set of all non-negative integers and $\beta_{j} =\frac{\sum_{i\in \mathcal{I}}\sum_{k\in \mathcal{K}}x_{ijk}}{\sum_{k\in \mathcal{K}}D_{jk}}$ is defined as the fill-rate of patient $j$ with $\beta$ being the maximum among $\beta_j$'s. Additionally, $u_i =\frac{\sum_{j\in \mathcal{J}}\sum_{k\in \mathcal{K}}x_{ijk}}{H_i}$ is the utilization level of caregiver $i$, and $u$ is the maximum among $u_i$'s. 

In Model \ref{eq:ObjFnc}, the objective function is represented by Equation \eqref{eq:ObjFnc}. The first term of the objective function is seeking to maximize the total utility (efficiency) of patients while penalizing the inequity defined by deviations from the maximum fill-rate $\theta(\beta-\beta_j)$ in the second term (achieving equity). The third term in the objective function is penalizing inefficacy by the deviations from the maximum utilization rate $\alpha(u-u_i)$ (achieving efficacy).

Constraint \eqref{eq:Const1} ensures that the total time allocated across the caregivers assigned to a patient for a care service does not exceed what the patient requires for that service. Constraint \eqref{eq:Const2} makes sure that a caregiver can serve a patient only if they are assigned to them. Constraint \eqref{eq:Const3} limits the total cost spent in the services less than the budget. Constraints \eqref{eq:Const4} and \eqref{eq:Const5} together ensure that a caregiver can be assigned to a patient for a service ($y_{ijk}=1$) only if they are capable of doing it ($e_{ik}=1$) and $z_{ij}$ equals 1 if caregiver $i$ provides at least one type of services to patient $j$ and 0 otherwise. Using the definition of the binary decision variable $z_{ij}$, constraint \eqref{eq:Const6} ensures that the total number of caregivers assigned to a patient cannot exceed the patient's preference level on the number of caregivers that take care of them. On the other hand, constraint \eqref{eq:Const7} makes sure that the total number of patients assigned to a caregiver does not exceed the caregiver's preference. Constraints \eqref{eq:Const8} and \eqref{eq:Const9} define fill-rate of each patient and utilization rate of each caregiver, respectively. Constraints \eqref{eq:Const10} and \eqref{eq:Const11} define the maximum satisfaction and maximum utilization rates, respectively. Finally, constraint \eqref{eq:Const12} sets the bounds and domains for decision variables.

\section{Model Analysis}
\label{sec:ModAnalysis}

In this section, we show that our problem defined by Model \ref{eq:ObjFnc} is NP-hard (non-deterministic polynomial-time hard)  by reducing the knapsack problem to our case. The knapsack problem (KP) can be formulated as an integer programming as follows.

\begin{align}
Max \hspace{3pt} \sum_{j\in\mathcal{J}} v_{j}x_{j} \label{eq:ObjFnc-KP}
\end{align}\setcounter{equation}{1}
\begin{subnumcases}{s.t.}
   \sum_{j\in\mathcal{J}}w_{j}x_{j} \leq W & \label{eq:Const1-KP}\\
   x_{j} \in \{0,1\}& $\forall j\in \mathcal{J}$, \label{eq:Const2-KP}
\end{subnumcases}
where, $\mathcal{J}$ is the set of items; $v_j$, $w_j$ and $W$ are item profit, item weight and  the knapsack capacity respectively; $x_{j}$ is decision variable ($x_{j}=1$, if item $j$ is selected; 0, otherwise).

It is well-known that the KP is NP-complete \citep{Kellerer2004knapsack}. If we can show that KP is polynomially reducible to our problem, then our problem is of course NP-hard. Otherwise, there exists polynomial-time algorithm for an NP-hard problem which is a contradiction. We conclude this in the following theorem and its proof. 

\begin{theorem}
\label{theo:NP-hard}
The caregivers assignment problem defined by Model \ref{eq:ObjFnc} is NP-hard. 
\end{theorem}

\noindent \underline{\emph{Proof}}

To prove the theorem, we show that the KP is polynomially reducible to our problem. To do this, we show any instance of the KP can be solved as a caregiver assignment problem. Given an instance
\begin{align}
\label{KP-concise}
    \max\{\sum_{j\in\mathcal{J}} v_{j}x_{j}: \; \sum_{j\in\mathcal{J}}w_{j}x_{j} \leq W, \, x_{j} \in \{0,1\}\},
\end{align}
we solve a caregiver assignment instance with only one service that is $\mathcal{I}=\mathcal{K}=\{1\}$; $\theta=\alpha=0$, $C=W$, and $e_1=1$, $N_j=D_j=1$, $M_1=H_1=|\mathcal{J}|$ for all $j$. Since there are only one caregiver and one type of services, the indices $i$ and $k$ are dropped. Furthermore, we let $p_j=v_j$ and $c_j=w_j$ for all $j$. 

Constraint \eqref{eq:Const1} reduces to $x_j \leq 1$ because $D_j=1$. However, considering $x_j$ is integer, this constraint now is equivalent to $x_j\in \{0, 1\}$. Because $C=W$ and $c_j=w_j$, constraint \eqref{eq:Const3} reduces to $\sum_{j\in \mathcal{J}} w_j x_j \leq W$. Constraint \eqref{eq:Const8} is equivalent to $\beta_j=x_j$ since $D_j=1$. Moreover, by the result that $x_j$ is binary and how we set the value of parameters, it is trivial to see that $x_j=y_j=z_j$ and all the other constraints are redundant. Finally, $\theta=\alpha=0$, $p_j=v_j$ and $\beta_j=x_j$ change the objective function to $\sum_{j\in\mathcal{J}} v_{j}x_{j}$.

Dropping the redundant constraints and variables, it leaves exactly the knapsack problem \eqref{KP-concise}. Therefore, an optimal solution of the caregiver assignment instance solves the knapsack instance. This implies that the KP is polynomially reducible to the caregiver assignment problem defined by Model \ref{eq:ObjFnc}.

\begin{flushright}   Q.E.D.   \end{flushright}

\begin{algorithm}[H]
\caption{Greedy Method (for building initial solution)}
\label{alg:greedy}
{\fontsize{10}{16}\selectfont
\begin{flushleft}
\textbf{Input: } 
\textbf{Output: } 
\end{flushleft}
	\begin{algorithmic}[1]
		\Procedure {Greedy}{$\mathcal{I}$, $\mathcal{J}$, $\mathcal{K}$, $C$, $M_i$'s, $N_j$'s, $H_i$'s, $\boldsymbol{e}_i$'s, $\boldsymbol{D}_j$'s, $c_{ijk}$'s, $p_j$'s, $\beta$}
		    \State sort $\mathcal{J}$ in descending order with key $p_j/\sum_{k\in \mathcal{K}} D_{jk}$
		    \State $c_0=C/\sum_{j\in\mathcal{J}} p_j$, $c_m = \displaystyle\min_{i\in \mathcal{I}, j\in \mathcal{J}, k\in \mathcal{K}} c_{ijk}$, $D_j=\sum_{k\in \mathcal{K}}D_{jk}$, $\mathcal{I}'=\emptyset$, $\mathscr{M}_i=\emptyset$, $\mathscr{N}_j=\emptyset$
			\For{$j \in \mathcal{J}$}
			    \State flag = False, $\boldsymbol{x}^*=\boldsymbol{0}$, $p^*=0$, $c^*=\infty$, $i^*=0$
			    \For{$i \in \mathcal{I}\backslash \mathcal{I}'$}
			        \State solve the optimization problem below to find $\boldsymbol{x}'=(x_1, \cdots, x_s)$ 
			        \begin{align}
			        \label{sp-greedy}
			            \max &\sum_{k\in \mathcal{K}} x'_{k}\\
			            s.t. &
			            \begin{cases}
			                           \sum_{k\in \mathcal{K}} c_{ijk}x'_{k}  \leq C &\\ x'_{k}\leq e_{ik}H_i, D_{jk} & \forall k \in \mathcal{K} \\
			                           \sum_{k\in \mathcal{K}} x'_{k} \leq H_i, \, \beta D_j &\\ x'_{k}\in \mathbb{N}_0 & \forall k\in \mathcal{K}
			            \end{cases} \notag
			        \end{align}
			        \State $p'=\frac{\sum_{k\in\mathcal{K}}x'_{k}}{\sum_{k\in\mathcal{K}}D_{jk}}p_j$, $c'=\sum_{k\in\mathcal{K}}c_{ijk} x'_{k}$
			        \If{$p'>p^*$ \textbf{and} $c'\leq p'c_0$ }
			            \State flag = True
			            \State $\boldsymbol{x}^*=\boldsymbol{x}'$, $p^*=p'$, $c^*=c'$, $i^*=i$
		            \ElsIf {flag = False \textbf{and} $p'>p^*$ \textbf{and} $c'\leq c^*$}
		                \State $\boldsymbol{x}^*=\boldsymbol{x}'$, $p^*=p'$, $c^*=c'$, $i^*=i$
	                \EndIf
                \EndFor
                
                \vspace{0.15cm}
                \State $M_{i^*}=M_{i^*}-1$, $N_j=N_j-1$, $C=C-c^*$, $\mathcal{I}'=\mathcal{I}'\cup \{i^*\}$, $D_j=D_j-\sum_{k\in \mathcal{K}}x^*_k/\beta$
                \State $\boldsymbol{x}_{i^*j}=\boldsymbol{x}^*$, $\boldsymbol{D}_j=\boldsymbol{D}_j-\boldsymbol{x}^*$, $p_j=p_j-p^*$, $H_{i^*}=H_{i^*}-\sum_{k\in\mathcal{K}}x^*_k$, $\mathscr{M}_{i^*}=\mathscr{M}_{i^*}\cup\{j\}$, $\mathscr{N}_{j}=\mathscr{N}_{j}\cup\{i^*\}$

                \If{$M_{i^*}=0$ \textbf{or} $H_{i^*}=0$}
                    \State $\mathcal{I}=\mathcal{I}\backslash\{i^*\}$
                \EndIf
                \If{$N_{j}=0$ \textbf{or} $\sum_{k\in\mathcal{K}} D_{jk}=0$}
                    \State $\mathcal{J}=\mathcal{J}\backslash\{j\}$
                \EndIf

    		    \If {$C \leq c_m$ \textbf{or} $|\mathcal{I}'|=n$} 
                    \State \textbf{break} 
                \EndIf
			\EndFor

	\algstore{myalg}
		
	\end{algorithmic}
	}
\end{algorithm}

\begin{algorithm}[H]
{\fontsize{10}{16}\selectfont
	\begin{algorithmic}
	\algrestore{myalg}
			\For{$j \in \mathcal{J}$}
			    \If{$C \leq c_m$ \textbf{or} $|\mathcal{I}|=0$}:
			        \State \textbf{break}
			    \EndIf
			    \For{$i \in \mathcal{I}$}
			        \State solve the sub-problem \eqref{sp-greedy} to find $\boldsymbol{x}'=(x_1, \cdots, x_s)$
			        \If{$\sum_{k \in \mathcal{K}} x'_k>0$ }
			            \State $c'=\sum_{k\in\mathcal{K}}c_{ijk} x'_{k}$
			            \State $M_{i}=M_{i}-1$, $N_j=N_j-1$, $C=C-c'$, $D_j=D_j-\sum_{k\in \mathcal{K}}x'_k/\beta$
			            \State $\boldsymbol{x}_{ij}=\boldsymbol{x}_{ij}+\boldsymbol{x}'$, $\boldsymbol{D}_j=\boldsymbol{D}_j-\boldsymbol{x}'$, $H_{i}=H_{i}-\sum_{k\in\mathcal{K}}x'_k$, $\mathscr{M}_{i}=\mathscr{M}_{i}\cup\{j\}$, $\mathscr{N}_{j}=\mathscr{N}_{j}\cup\{i\}$
                        \If{$M_{i}=0$ \textbf{or} $H_{i}=0$}
                            \State $\mathcal{I}=\mathcal{I}\backslash\{i\}$
                        \EndIf
                        \If{$N_{j}=0$ \textbf{or} $\sum_{k\in\mathcal{K}} D_{jk}=0$ \textbf{or} $C \leq c_m$}
                            \State $\mathcal{J}=\mathcal{J}\backslash\{j\}$
                            \State \textbf{break}
                        \EndIf
	                \EndIf
                \EndFor
			\EndFor
			
			\Return $\boldsymbol{x}_{ij}$'s, $H_{i}$'s, $\boldsymbol{D}_j$'s, $\mathscr{M}_{i}$'s, $\mathscr{N}_{j}$'s
		\EndProcedure

	\end{algorithmic}
	}
\end{algorithm}

\begin{algorithm}[H]
\caption{Tabu Search (for finding the local optima)}
\label{alg:tabu}
{\fontsize{10}{16}\selectfont
\begin{flushleft}
\textbf{Input: } the outputs from Algorithm  \ref{alg:greedy} and the parameters $\mathcal{I}$, $\mathcal{J}$, $\mathcal{K}$, $\boldsymbol{e}_i$'s, $\beta$ \\
\textbf{Output: } \\
\textbf{Initialization: } set tabu lists $\mathcal{T}_c=\{ \}$ and $\mathcal{T}_p=\{ \}$, the current time $t=0$, and the maximum running time length to be $T$.
\end{flushleft}
	\begin{algorithmic}[1]
		\Procedure {Tabu}{outputs from Algorithm  \ref{alg:greedy}, $\mathcal{I}$, $\mathcal{J}$, $\mathcal{K}$, $\boldsymbol{e}_i$'s, $\beta$}
		    \State calculate $\beta_j$ for any $j \in \mathcal{J}$ and let $\mathcal{J}'=\{j \mid \beta_j<\beta, \, j\in \mathcal{J}\}$
		    \State sort $\mathcal{J}'$ in ascending order with key $\beta_j$
		    \While {$t<T$}
		        \While {$\mathcal{J}'\neq \emptyset$}
		        	\State $\mathcal{I}'=\{i\mid H_i>0, \, i \in \mathcal{I}\}$
		            \State sort $\mathcal{I}'$ in descending order with key $\min \{H_i, \sum_{k\in \mathcal{K}} e_{ik}D_{jk}\}$
		            \For{$j \in \mathcal{J}'$}
    		            \If {$\mathcal{I}'=\emptyset$}
    		                \State $\mathcal{J}'=\emptyset$
    		                \State \textbf{break}
    		            \EndIf

		            \EndFor

		        \EndWhile
		        \If {$\mathcal{J}'=\emptyset$}
		            \State \textbf{break}
		        \EndIf
		    \EndWhile

	\algstore{myalg1}
		
	\end{algorithmic}
	}
\end{algorithm}

\begin{algorithm}[H]
{\fontsize{10}{16}\selectfont
	\begin{algorithmic}
	\algrestore{myalg1}
			\For{$j \in \mathcal{J}$}
			    \If{$C \leq c_m$ \textbf{or} $|\mathcal{I}|=0$}:
			        \State \textbf{break}
			    \EndIf
			    \For{$i \in \mathcal{I}$}
			        \State solve the sub-problem \eqref{sp-greedy} to find $\boldsymbol{x}'=(x_1, \cdots, x_s)$
			        \If{$\sum_{k \in \mathcal{K}} x'_k>0$ }
			            \State $c'=\sum_{k\in\mathcal{K}}c_{ijk} x'_{k}$
			            \State $M_{i}=M_{i}-1$, $N_j=N_j-1$, $C=C-c'$, $D_j=D_j-\sum_{k\in \mathcal{K}}x'_k/\beta$
			            \State $\boldsymbol{x}_{ij}=\boldsymbol{x}_{ij}+\boldsymbol{x}'$, $\boldsymbol{D}_j=\boldsymbol{D}_j-\boldsymbol{x}'$, $H_{i}=H_{i}-\sum_{k\in\mathcal{K}}x'_k$
                        \If{$M_{i}=0$ \textbf{or} $H_{i}=0$}
                            \State $\mathcal{I}=\mathcal{I}\backslash\{i\}$
                        \EndIf
                        \If{$N_{j}=0$ \textbf{or} $\sum_{k\in\mathcal{K}} D_{jk}=0$ \textbf{or} $C \leq c_m$}
                            \State $\mathcal{J}=\mathcal{J}\backslash\{j\}$
                            \State \textbf{break}
                        \EndIf
	                \EndIf
                \EndFor
			\EndFor
			
			\Return $\boldsymbol{x}_{ij}$'s
		\EndProcedure

	\end{algorithmic}
	}
\end{algorithm}

\section{Numerical Study}
\label{sec:NumStudy}

\section{Conclusions}
\label{sec:Conclusions}

\newpage
\begin{appendix}

\end{appendix}

\newpage
\bibliographystyle{ormsv080}
\bibliography{main}

\end{document}